\documentclass[twoside,reqno]{amsart}

\usepackage{latexsym}

\newcommand{\qdn}{\hspace*{-1.5mm}}
\newcommand{\qqdn}{\hspace*{-2.5mm}}
\newcommand{\xqdn}{\hspace*{-5.0mm}}
\newcommand{\xxqdn}{\hspace*{-10mm}}

\newcommand{\fns}{\footnotesize}

\newcommand{\fnk}[3]{\left[\qdn\ba{#1}#2\\#3\ea\qdn\right]}
\newcommand{\ffnk}[4]{\left[\qdn\ba{#1}#3\\#4\ea{\!\Big|\:#2}\right]}

\newcommand{\binm}{\binom}
\newcommand{\sbnm}[2]{\Bigl(\!\ba{c}\!#1\!\\#2\ea\!\Bigr)}

\newcommand{\nnm}{\nonumber}
\newcommand{\be}{\begin{equation}}
\newcommand{\ee}{\end{equation}}
\newcommand{\ba}{\begin{array}}
\newcommand{\ea}{\end{array}}
\newcommand{\bmn}{\begin{eqnarray}}
\newcommand{\emn}{\end{eqnarray}}
\newcommand{\bnm}{\begin{eqnarray*}}
\newcommand{\enm}{\end{eqnarray*}}
\newcommand{\bln}{\begin{subequations}}
\newcommand{\eln}{\end{subequations}}

\newtheorem{thm}{Theorem}

\newtheorem{corl}[thm]{Corollary}
\newtheorem{prop}[thm]{Proposition}
\newtheorem{exam}{Example}

\newtheorem{entry}{Entry}

\newcommand{\bbtm}[4]{\bibitem{kn:#1}{#2,}~{#3,}~{#4.}}
\newcommand{\cito}[1]{\cite{kn:#1}}
\newcommand{\citu}[2]{\cite[#2]{kn:#1}}

\begin{document} 
{\fns
\title{Derivative operator \\and harmonic number identities}
\author{$^A$Chuanan Wei and $^B$Dianxuan Gong }
\dedicatory{$^A$Department of Information Technology\\
            Hainan Medical College, Haikou 571101, China\\
            $^B$College of Sciences\\
             Hebei Polytechnic University, Tangshan 063009, China }
\thanks{\emph{Email addresses}: weichuanan@yahoo.com.cn (C. Wei),
        dxgong@heut.edu.cn (D. Gong)}

\address{ }


\begin{abstract}
 By applying the derivative operator to the corresponding hypergeometric form of a $q$-series
 transformation due to Andrews \citu{andrews-a}{Theorem 4},
  we establish a general harmonic number identity.  As the special
cases of it, several interesting Chu-Donno type identities and
Paule-Schneider type identities are displayed.
\end{abstract}

\maketitle\thispagestyle{empty}
\markboth{C. Wei and D. Gong}
         {Derivative operator and harmonic number identities}

\section{Introduction}



For a nonnegative integer $n$, define the harmonic numbers by
\[H_{0}=0\quad \text{and}\quad H_{n}=\sum_{k=1}^n\frac{1}{k}\quad \text{when}\quad n=1,2,\cdots.\]
For a differentiable function $f(x)$, the derivative
 operator $\mathcal{D}$  can be defined by
\[\mathcal{D}f(x)=\frac{d}{dx}f(x)\Big|_{x=0}.\]
Then it is not difficult to show the  following two derivatives of
binomial coefficients:
 \bnm
&&\:\mathcal{D}\:\binm{n+x}{r}=\binm{n}{r}\big\{H_n-H_{n-r}\big\},\\
&&\:\mathcal{D}\:\binm{n-x}{r}=\binm{n}{r}\big\{H_{n-r}-H_{n}\big\},
 \enm
where $r\leq n$ with $r=0,1,\cdots$.

 For a complex number $x$, define the shifted factorial by
  \[(x)_0=1\quad\text{and}\quad (x)_n=\prod_{k=0}^{n-1}(x+k)
  \quad\text{when}\quad n=1,2,\cdots.\]
The fractional form of it reads as
\[\fnk{ccccc}{a,&b,&\cdots,&c}{\alpha,&\beta,&\cdots,&\gamma}_n
=\frac{(a)_n(b)_n\cdots(c)_n}{(\alpha)_n(\beta)_n\cdots(\gamma)_n}.\]
Then the hypergeometric series(cf. Bailey~\cito{bailey}) can be
defined by
\[_{1+r}F_s\ffnk{cccc}{z}{a_0,&a_1,&\cdots,&a_r}
{&b_1,&\cdots,&b_s} \:=\:\sum_{k=0}^\infty
\fnk{ccccc}{a_0,&a_1,&\cdots,&a_r}{1,&b_1,&\cdots,&b_s}_kz^k,\]
where $\{a_{i}\}_{i\geq0}$ and $\{b_{j}\}_{j\geq1}$ are complex
parameters such that no zero factors appear in the denominators of
the summand on the right hand side.

 For a complex sequence $\{A_k\}_{k\geq0}$ and two nonnegative integers $i$ and $j$,
define the product by
 \bnm
\prod_{k=i}^jA_k=
\begin{cases}
A_iA_{i+1}\cdots A_j,\quad\:\: \text{for}\quad j\geq i,\\
\qquad\quad1,\qquad\quad\quad \text{for}\quad j=i-1.
\end{cases}
 \enm

In 1975, Andrews \citu{andrews-a}{Theorem 4} gave a beautiful
$q$-series transformation. The corresponding hypergeometric form of
it(cf. Krattenthaler et al. \citu{krattenthaler-a}{Theorem 8} and
\citu{krattenthaler-c}{Equation 4.2}) can be stated as
 \bmn
&&\xqdn\qdn_{2m+5}F_{2m+4}\ffnk{cccc}{1}{a,&\qdn1+a/2,&\qdn\{P_s\}_{s=1}^{2m+2},&\qdn-n}
{&\qdn a/2,&\qdn\{1+a-P_{s}\}_{s=1}^{2m+2},&\qdn1+a+n}=
\fnk{ccccc}{1+a,\:1+a-P_{2m+1}-P_{2m+2}}{1+a-P_{2m+1},\:1+a-P_{2m+2}}_n
 \nnm\\\label{source}
 &&\times\qdn\sum_{0\leq i_1\leq i_2\leq\cdots\leq i_m \leq n}
\prod_{r=1}^{m}\fnk{ccccc}{-i_{r+1},\quad P_{2r+1},\quad
P_{2r+2},\quad1+a-P_{2r-1}-P_{2r}} {1,\quad
P_{2r+1}+P_{2r+2}-a-i_{r+1},\quad1+a-P_{2r-1},\quad1+a-P_{2r}}_{i_{r}},
\emn
  where $i_{m+1}=n$ and $m\in\mathbb{N}$. When $m=1$, the last equation reduces to the famous Whipple's transformation
(cf. Bailey \citu{bailey}{p. 25}):
 \bnm
 &&\xxqdn\qdn{_7F_6}\ffnk{ccccccccc}{1}
 {a,&1+a/2,&P_1,&P_2,&P_3,&P_4,&-n}
 {&a/2,&1+a-P_1,&1+a-P_2,&1+a-P_3,&1+a-P_4,&1+a+n}
  \nnm\\ \label{whipple}
&&\xxqdn=\:\fnk{ccc}{1+a,1+a-P_3-P_4}{1+a-P_3,1+a-P_4}_n
{_4F_3}\fnk{cccc}{-n,\quad P_3,\quad P_4,\quad 1+a-P_1-P_2}
{P_3+P_4-a-n,\quad 1+a-P_1,\quad1+a-P_2}.
 \enm

 By applying the derivative operator $D$ to \eqref{source}, we shall establish a general
harmonic number identity in the next section.  As the special cases
of it, several interesting Chu-Donno type identities and
Paule-Schneider type identities will  be displayed.

\section{Harmonic number identities}
\subsection{A general harmonic number identity}${}$

Let $v$ be a nonnegative integer with $0\leq v\leq 2m+2$. For two
finite sets $\{\alpha_s\}_{s=1}^{v}$ and
$\{\alpha_s\}_{s=\,v+1}^{2m+2}$, the case $v=0$ corresponds to the
former is empty and the latter is $\{\alpha_s\}_{s=1}^{2m+2}$, and
the case $v=2m+2$ corresponds to the former is
$\{\alpha_s\}_{s=1}^{2m+2}$ and the latter is empty.

Performing the replacements $a\to-x-n$, $P_s \to 1+P_s$ with $1\leq
s\leq v$ and $P_s \to -n-P_s$ with $v+1\leq s\leq 2m+2$ for
\eqref{source}, we obtain the following expression:
 \bnm
&&\xqdn\qqdn\sum_{k=0}^n\binm{n}{k}(x+n-2k)\frac{\binm{n+x}{k}}{\binm{k-x}{k}}
\prod_{s=1}^{v}\frac{\binm{k+P_s}{k}}{\binm{x+n+P_s}{k}}
\prod_{s=\,v+1}^{2m+2}\frac{\binm{n+P_s}{k}}{\binm{-x+k+P_s}{k}}\\
&&\xqdn\qdn=\:\:x\:\fnk{ccccc}{-x-n,\quad1-x-n-T_{2m+1}-T_{2m+2}}{1-x-n-T_{2m+1},\quad1-x-n-T_{2m+2}}_n
 \\&&\qdn\xqdn\times\qqdn
\sum_{0\leq i_1\leq i_2\leq\cdots\leq i_m \leq n} \prod_{r=1}^{m}
\fnk{ccccc}{-i_{r+1},\:\:T_{2r+1},\:\:T_{2r+2},\:\:1-x-n-T_{2r-1}-T_{2r}}
{1,\:\:T_{2r+1}+T_{2r+2}+x+n-i_{r+1},\:\:1-x-n-T_{2r-1},\:\:1-x-n-T_{2r}}_{i_r},
 \enm
 where $T_s=1+P_s$ with $1\leq s\leq v$ and $T_s=-n-P_s$ with $v+1\leq s\leq
 2m+2$.

Applying the derivative operator $\mathcal{D}$ to both sides of the
last equation, we establish the following theorem.

\begin{thm} \label{thm}
For $2m+2$ nonnegative integers $\{P_s\}_{s=1}^{2m+2}$  with
$i_{m+1}=n$ and $m\in\mathbb{N}$, there holds the general harmonic
number identity:
 \bnm
&&\sum_{k=0}^n\binm{n}{k}^2
 \prod_{s=1}^{v}\frac{\binm{k+P_s}{k}}{\binm{n+P_s}{k}}
\prod_{s=\,v+1}^{2m+2}\frac{\binm{n+P_s}{k}}{\binm{k+P_s}{k}}
\Bigg\{1+(n-2k)\bigg(2H_k-\sum_{s=1}^{v}H_{k+P_s}+\sum_{s=v+1}^{2m+2}H_{k+P_s}\bigg)\Bigg\}
\\&&\:\:=\:\fnk{ccccc}{-n,\quad1-n-T_{2m+1}-T_{2m+2}}{1-n-T_{2m+1},\quad1-n-T_{2m+2}}_n
 \\&&\:\:\times\:
\sum_{0\leq i_1\leq i_2\leq\cdots\leq i_m \leq n} \prod_{r=1}^{m}
\fnk{ccccc}{-i_{r+1},\:\:T_{2r+1},\:\:T_{2r+2},\:\:1-n-T_{2r-1}-T_{2r}}
{1,\:\:T_{2r+1}+T_{2r+2}+n-i_{r+1},\:\:1-n-T_{2r-1},\:\:1-n-T_{2r}}_{i_r},
 \enm
 where $T_s=1+P_s$ with $1\leq s\leq v$ and $T_s=-n-P_s$ with $v+1\leq s\leq
 2m+2$.
\end{thm}

\subsection{Special cases: harmonic number identities of Chu-Donno
type}${}$

 Setting $m=1$ in Theorem \ref{thm}, we get the
following equation.

\begin{prop}\label{prop-a}
For four nonnegative integers $\{P_s\}_{s=1}^4$,  there holds the
harmonic number identity:
 \bnm
 &&\xxqdn\sum_{k=0}^n\binm{n}{k}^2
 \prod_{s=1}^{v}\frac{\binm{k+P_s}{k}}{\binm{n+P_s}{k}}
\prod_{s=\,v+1}^{4}\frac{\binm{n+P_s}{k}}{\binm{k+P_s}{k}}
\Bigg\{1+(n-2k)\bigg(2H_k-\sum_{s=1}^{v}H_{k+P_s}+\sum_{s=v+1}^{4}H_{k+P_s}\bigg)\Bigg\}
\\&&\xxqdn\:\:=\:\fnk{ccccc}{-n,\quad1-n-T_{3}-T_{4}}{1-n-T_{3},\quad1-n-T_{4}}_n
{_{4}F_3}\ffnk{ccccc}{1}{-n,\quad T_{3},\quad
T_{4},\quad1-n-T_{1}-T_{2}}{\quad
T_{3}+T_{4},\quad1-n-T_{1},\quad1-n-T_{2}},
 \enm
 where $T_s=1+P_s$ with $1\leq s\leq v$ and $T_s=-n-P_s$ with $v+1\leq s\leq
 4$.
\end{prop}

Letting $P_s\to nP_s$ with $1\leq s\leq4$ in Proposition
\ref{prop-a}, we get the following result.

\begin{corl}\label{corl}
For four nonnegative integers $\{P_s\}_{s=1}^4$, there holds the
harmonic number identity:
 \bnm
 &&\xxqdn\sum_{k=0}^n\binm{n}{k}^2
 \prod_{s=1}^{v}\frac{\binm{k+nP_s}{k}}{\binm{n+nP_s}{k}}
\prod_{s=\,v+1}^{4}\frac{\binm{n+nP_s}{k}}{\binm{k+nP_s}{k}}
\Bigg\{1+(n-2k)\bigg(2H_k-\sum_{s=1}^{v}H_{k+nP_s}+\sum_{s=v+1}^{4}H_{k+nP_s}\bigg)\Bigg\}
\\&&\xxqdn\:\:=\:\fnk{ccccc}{-n,\quad1-n-T_{3}-T_{4}}{1-n-T_{3},\quad1-n-T_{4}}_n
{_{4}F_3}\ffnk{ccccc}{1}{-n,\quad T_{3},\quad
T_{4},\quad1-n-T_{1}-T_{2}}{\quad
T_{3}+T_{4},\quad1-n-T_{1},\quad1-n-T_{2}},
 \enm
 where $T_s=1+nP_s$ with $1\leq s\leq v$ and $T_s=-n-nP_s$ with $v+1\leq s\leq
 4$.
\end{corl}

The importance of Corollary \ref{corl} lies in that it implies eight
important theorems due to Chu and Donno \cito{chu}. The details are
laid out as follows.

Taking  respectively $v=\:2,\:1,\:0$ in Corollary \ref{corl} and
then letting $P_1\to b$, $P_2\to c$, $P_3\to \infty$, $P_4\to
\infty$, we gain the following three known harmonic number
identities.

\begin{exam}[Chu and Donno \citu{chu}{Theorem 5}]\label{exam-a}
 For two nonnegative integers $\{b,c\}$, there holds
\bnm
 \xqdn\qqdn\sum_{k=0}^n\binm{n}{k}^2
 \frac{\binm{k+bn}{k}\binm{k+cn}{k}}
{\binm{n+bn}{k}\binm{n+cn}{k}}
\Big\{1+(n-2k)(2H_k-H_{bn+k}-H_{cn+k})\Big\}
=\frac{\binm{1+bn+cn+n}{n}}{\binm{n+bn}{n}\binm{n+cn}{n}}.
  \enm
\end{exam}

\begin{exam}[Chu and Donno \citu{chu}{Theorem 6}]\label{exam-b}
 For two nonnegative integers $\{b,c\}$,there holds
\bnm
 \sum_{k=0}^n\binm{n}{k}^2
 \frac{\binm{k+bn}{k}\binm{n+cn}{k}}
{\binm{n+bn}{k}\binm{k+cn}{k}}
\Big\{1+(n-2k)(2H_k-H_{bn+k}+H_{cn+k})\Big\}
=(-1)^n\frac{\binm{bn-cn}{n}}{\binm{n+bn}{n}\binm{n+cn}{n}}.
  \enm
\end{exam}

\begin{exam}[Chu and Donno \citu{chu}{Theorem 7}]\label{exam-c}
 For two nonnegative integers $\{b,c\}$, there holds
\bnm
 \sum_{k=0}^n\binm{n}{k}^2
 \frac{\binm{n+bn}{k}\binm{n+cn}{k}}
{\binm{k+bn}{k}\binm{k+cn}{k}}
\Big\{1+(n-2k)(2H_k+H_{bn+k}+H_{cn+k})\Big\}
=(-1)^n\frac{\binm{2n+bn+cn}{n}}{\binm{n+bn}{n}\binm{n+cn}{n}}.
  \enm
\end{exam}

Taking respectively $v=4,\:3,\:2,\:1,\:0$ in Corollary \ref{corl}
and then letting $P_1\to b$, $P_2\to c$, $P_3\to d$, $P_4\to e\:$,
we achieve the following five known harmonic number identities.

\begin{exam}[Chu and Donno \citu{chu}{Theorem 8}]\label{exam-d}
 For four nonnegative integers $\{b,c,d,e\}$, there holds
\bnm
 && \xxqdn\xqdn\sum_{k=0}^n\binm{n}{k}^2
 \frac{\binm{k+bn}{k}\binm{k+cn}{k}\binm{k+dn}{k}\binm{k+en}{k}}
{\binm{n+bn}{k}\binm{n+cn}{k}\binm{n+dn}{k}\binm{n+en}{k}}
\\&&\xxqdn\xqdn\:\:\times\:
\Big\{1+(n-2k)(2H_k-H_{bn+k}-H_{cn+k}-H_{dn+k}-H_{en+k})\Big\}\\
&&\xxqdn\xqdn\:\:=\:
 \frac{\binm{1+bn+cn+n}{n}}{\binm{n+bn}{n}\binm{n+cn}{n}}
\sum_{i=0}^n\binm{n}{i}\frac{\binm{i+bn}{i}\binm{i+cn}{i}\binm{1+dn+en+n}{i}}
{\binm{n+dn}{i}\binm{n+en}{i}\binm{1+bn+cn+i}{i}}.
  \enm
\end{exam}

\begin{exam}[Chu and Donno \citu{chu}{Theorem 9}]\label{exam-e}
 For four nonnegative integers $\{b,c,d,e\}$, there holds
\bnm
 && \xxqdn\xqdn\sum_{k=0}^n\binm{n}{k}^2
 \frac{\binm{k+bn}{k}\binm{k+cn}{k}\binm{k+dn}{k}\binm{n+en}{k}}
{\binm{n+bn}{k}\binm{n+cn}{k}\binm{n+dn}{k}\binm{k+en}{k}}
\\&&\xxqdn\xqdn\:\:\times\:
\Big\{1+(n-2k)(2H_k-H_{bn+k}-H_{cn+k}-H_{dn+k}+H_{en+k})\Big\}
\\&&\xxqdn\xqdn\:\:=\:\frac{\binm{1+bn+cn+n}{n}}{\binm{n+bn}{n}\binm{n+cn}{n}}
\sum_{i=0}^n(-1)^i\binm{n}{i}\frac{\binm{i+bn}{i}\binm{i+cn}{i}\binm{dn-en}{i}}
{\binm{n+dn}{i}\binm{i+en}{i}\binm{1+bn+cn+i}{i}}.
  \enm
\end{exam}

\begin{exam}[Chu and Donno \citu{chu}{Theorem 10}]\label{exam-f}
 For four nonnegative integers $\{b,c,d,e\}$, there holds
\bnm
 &&\xxqdn\xqdn \sum_{k=0}^n\binm{n}{k}^2
 \frac{\binm{k+bn}{k}\binm{k+cn}{k}\binm{n+dn}{k}\binm{n+en}{k}}
{\binm{n+bn}{k}\binm{n+cn}{k}\binm{k+dn}{k}\binm{k+en}{k}}
\\&&\xxqdn\xqdn\:\:\times\:
\Big\{1+(n-2k)(2H_k-H_{bn+k}-H_{cn+k}+H_{dn+k}+H_{en+k})\Big\}\\
&&\xxqdn\xqdn\:\:=\:
\frac{\binm{1+bn+cn+n}{n}}{\binm{n+bn}{n}\binm{n+cn}{n}}
\sum_{i=0}^n(-1)^i\binm{n}{i}\frac{\binm{i+bn}{i}\binm{i+cn}{i}\binm{n+dn+en+i}{i}}
{\binm{i+dn}{i}\binm{i+en}{i}\binm{1+bn+cn+i}{i}}.
  \enm
\end{exam}

\begin{exam}[Chu and Donno \citu{chu}{Theorem 11}]\label{exam-g}
 For four nonnegative integers $\{b,c,d,e\}$, there holds
\bnm
 &&\xxqdn\xqdn \sum_{k=0}^n\binm{n}{k}^2
 \frac{\binm{k+bn}{k}\binm{n+cn}{k}\binm{n+dn}{k}\binm{n+en}{k}}
{\binm{n+bn}{k}\binm{k+cn}{k}\binm{k+dn}{k}\binm{k+en}{k}}
 \\&&\xxqdn\xqdn\:\:\times\:
\Big\{1+(n-2k)(2H_k-H_{bn+k}+H_{cn+k}+H_{dn+k}+H_{en+k})\Big\}
 \\&&\xxqdn\xqdn\:\:=\:
(-1)^n\frac{\binm{bn-cn}{n}}{\binm{n+bn}{n}\binm{n+cn}{n}}
\sum_{i=0}^n\binm{n}{i}\frac{\binm{i+bn}{i}\binm{n+cn}{i}\binm{i+dn+en+n}{i}}
{\binm{i+dn}{i}\binm{i+en}{i}\binm{i+bn-cn-n}{i}}.
  \enm
\end{exam}

\begin{exam}[Chu and Donno \citu{chu}{Theorem 12}]\label{exam-h}
 For four nonnegative integers $\{b,c,d,e\}$, there holds
\bnm
 &&\xxqdn\xqdn \sum_{k=0}^n\binm{n}{k}^2
 \frac{\binm{n+bn}{k}\binm{n+cn}{k}\binm{n+dn}{k}\binm{n+en}{k}}
{\binm{k+bn}{k}\binm{k+cn}{k}\binm{k+dn}{k}\binm{k+en}{k}}
 \\&&\xxqdn\xqdn\:\:\times\:
\Big\{1+(n-2k)(2H_k+H_{bn+k}+H_{cn+k}+H_{dn+k}+H_{en+k})\Big\}\\
&&\xxqdn\xqdn\:\:=\:
(-1)^n\frac{\binm{2n+bn+cn}{n}}{\binm{n+bn}{n}\binm{n+cn}{n}}
\sum_{i=0}^n\binm{n}{i}\frac{\binm{n+bn}{i}\binm{n+cn}{i}\binm{n+dn+en+i}{i}}
{\binm{i+dn}{i}\binm{i+en}{i}\binm{2n+bn+cn}{i}}.
  \enm
\end{exam}

Setting $m=2$ in Theorem \ref{thm}, we attain the following
equation.

\begin{prop}\label{prop-b}
For six nonnegative integers $\{P_s\}_{s=1}^6$, there holds the
harmonic number identity:
 \bnm
 &&\xxqdn\sum_{k=0}^n\binm{n}{k}^2
 \prod_{s=1}^{v}\frac{\binm{k+P_s}{k}}{\binm{n+P_s}{k}}
\prod_{s=\,v+1}^{6}\frac{\binm{n+P_s}{k}}{\binm{k+P_s}{k}}
\Bigg\{1+(n-2k)\bigg(2H_k-\sum_{s=1}^{v}H_{k+P_s}+\sum_{s=v+1}^{6}H_{k+P_s}\bigg)\Bigg\}
\\&&\xxqdn\:\:=\:
 \fnk{ccccc}{-n,\quad1-n-T_{5}-T_{6}}{1-n-T_{5},\quad1-n-T_{6}}_n
\sum_{i=0}^n \fnk{ccccc}{-n,\quad T_{5},\quad
T_{6},\quad1-n-T_{3}-T_{4}} {1,\quad
T_{5}+T_{6},\quad1-n-T_{3},\quad1-n-T_{4}}_{i}
 \\&&\xxqdn\:\:\:\times\:\:
 {_{4}F_3}\ffnk{ccccc}{1}{-i,\quad
T_{3},\quad T_{4},\quad1-n-T_{1}-T_{2}}{\quad
T_{3}+T_{4}+n-i,\quad1-n-T_{1},\quad1-n-T_{2} },
 \enm
 where $T_s=1+P_s$ with $1\leq s\leq v$ and $T_s=-n-P_s$ with $v+1\leq s\leq
 6$.
\end{prop}

Taking respectively $v=6,\:5,\:4,\:3,\:2,\:1,\:0$ in Proposition
\ref{prop-b} and then letting $P_1\to b$, $P_2\to c$, $P_3\to d$,
$P_4\to e$ $P_5\to f$, $P_6\to g$, we derive the following seven
harmonic number identities  of Chu-Donno type.

\begin{exam}\label{exam-i}
 For six nonnegative integers $\{b,c,d,e,f,g\}$, there holds
\bnm
 &&\xxqdn\qdn \sum_{k=0}^n\binm{n}{k}^2
 \frac{\binm{k+b}{k}\binm{k+c}{k}\binm{k+d}{k}\binm{k+e}{k}\binm{k+f}{k}\binm{k+g}{k}}
{\binm{n+b}{k}\binm{n+c}{k}\binm{n+d}{k}\binm{n+e}{k}\binm{n+f}{k}\binm{n+g}{k}}
\\&&\xxqdn\qdn\:\:\times\:
\Big\{1+(n-2k)(2H_k-H_{b+k}-H_{c+k}-H_{d+k}-H_{e+k}-H_{f+k}-H_{g+k})\Big\}\\
&&\xxqdn\qdn\:\:=\:\frac{\binm{1+b+c+n}{n}}{\binm{n+b}{n}\binm{n+c}{n}}
\sum_{i=0}^n\binm{n}{i}\frac{\binm{i+b}{i}\binm{i+c}{i}\binm{1+d+e+n}{i}}
{\binm{n+d}{i}\binm{n+e}{i}\binm{1+b+c+i}{i}}
\sum_{j=0}^i\binm{i}{j}\frac{\binm{j+d}{j}\binm{j+e}{j}\binm{1+f+g+n}{j}}
{\binm{n+f}{j}\binm{n+g}{j}\binm{1+d+e+n-i+j}{j}}.
  \enm
\end{exam}

\begin{exam}\label{exam-j}
 For six nonnegative integers $\{b,c,d,e,f,g\}$,  there holds
\bnm
 &&\xqdn \sum_{k=0}^n\binm{n}{k}^2
 \frac{\binm{k+b}{k}\binm{k+c}{k}\binm{k+d}{k}\binm{k+e}{k}\binm{k+f}{k}\binm{n+g}{k}}
{\binm{n+b}{k}\binm{n+c}{k}\binm{n+d}{k}\binm{n+e}{k}\binm{n+f}{k}\binm{k+g}{k}}
\\&&\xqdn\:\:\times\:
\Big\{1+(n-2k)(2H_k-H_{b+k}-H_{c+k}-H_{d+k}-H_{e+k}-H_{f+k}+H_{g+k})\Big\}\\
&&\xqdn\:\:=\:\frac{\binm{1+b+c+n}{n}}{\binm{n+b}{n}\binm{n+c}{n}}
\sum_{i=0}^n\binm{n}{i}\frac{\binm{i+b}{i}\binm{i+c}{i}\binm{1+d+e+n}{i}}
{\binm{n+d}{i}\binm{n+e}{i}\binm{1+b+c+i}{i}}
\sum_{j=0}^i(-1)^j\binm{i}{j}\frac{\binm{j+d}{j}\binm{j+e}{j}\binm{f-g}{j}}
{\binm{n+f}{j}\binm{j+g}{j}\binm{1+d+e+n-i+j}{j}}.
  \enm
\end{exam}

\begin{exam}\label{exam-k}
 For six nonnegative integers $\{b,c,d,e,f,g\}$,  there holds
\bnm
 &&\xqdn \sum_{k=0}^n\binm{n}{k}^2
 \frac{\binm{k+b}{k}\binm{k+c}{k}\binm{k+d}{k}\binm{k+e}{k}\binm{n+f}{k}\binm{n+g}{k}}
{\binm{n+b}{k}\binm{n+c}{k}\binm{n+d}{k}\binm{n+e}{k}\binm{k+f}{k}\binm{k+g}{k}}
\\&&\xqdn\:\:\times\:
\Big\{1+(n-2k)(2H_k-H_{b+k}-H_{c+k}-H_{d+k}-H_{e+k}+H_{f+k}+H_{g+k})\Big\}\\
&&\xqdn\:\:=\:\frac{\binm{1+b+c+n}{n}}{\binm{n+b}{n}\binm{n+c}{n}}
\sum_{i=0}^n\binm{n}{i}\frac{\binm{i+b}{i}\binm{i+c}{i}\binm{1+d+e+n}{i}}
{\binm{n+d}{i}\binm{n+e}{i}\binm{1+b+c+i}{i}}
\sum_{j=0}^i(-1)^j\binm{i}{j}\frac{\binm{j+d}{j}\binm{j+e}{j}\binm{n+f+g+j}{j}}
{\binm{j+f}{j}\binm{j+g}{j}\binm{1+d+e+n-i+j}{j}}.
  \enm
\end{exam}

\begin{exam}\label{exam-l}
 For six nonnegative integers $\{b,c,d,e,f,g\}$,  there holds
\bnm
 &&\xqdn\qqdn \sum_{k=0}^n\binm{n}{k}^2
 \frac{\binm{k+b}{k}\binm{k+c}{k}\binm{k+d}{k}\binm{n+e}{k}\binm{n+f}{k}\binm{n+g}{k}}
{\binm{n+b}{k}\binm{n+c}{k}\binm{n+d}{k}\binm{k+e}{k}\binm{k+f}{k}\binm{k+g}{k}}
\\&&\xqdn\qqdn\:\:\times\:
\Big\{1+(n-2k)(2H_k-H_{b+k}-H_{c+k}-H_{d+k}+H_{e+k}+H_{f+k}+H_{g+k})\Big\}\\
&&\xqdn\qqdn\:\:=\:\frac{\binm{1+b+c+n}{n}}{\binm{n+b}{n}\binm{n+c}{n}}
\sum_{i=0}^n(-1)^i\binm{n}{i}\frac{\binm{i+b}{i}\binm{i+c}{i}\binm{d-e}{i}}
{\binm{n+d}{i}\binm{i+e}{i}\binm{1+b+c+i}{i}}
\sum_{j=0}^i\binm{i}{j}\frac{\binm{j+d}{j}\binm{n+e}{j}\binm{j+f+g+n}{j}}
{\binm{j+f}{j}\binm{j+g}{j}\binm{j+d-e-i}{j}}.
  \enm
\end{exam}

\begin{exam}\label{exam-m}
 For six nonnegative integers $\{b,c,d,e,f,g\}$, there holds
\bnm
 &&\xqdn\qqdn \sum_{k=0}^n\binm{n}{k}^2
 \frac{\binm{k+b}{k}\binm{k+c}{k}\binm{n+d}{k}\binm{n+e}{k}\binm{n+f}{k}\binm{n+g}{k}}
{\binm{n+b}{k}\binm{n+c}{k}\binm{k+d}{k}\binm{k+e}{k}\binm{k+f}{k}\binm{k+g}{k}}
\\&&\xqdn\qqdn\:\:\times\:
\Big\{1+(n-2k)(2H_k-H_{b+k}-H_{c+k}+H_{d+k}+H_{e+k}+H_{f+k}+H_{g+k})\Big\}\\
&&\xqdn\qqdn\:\:=\:\frac{\binm{1+b+c+n}{n}}{\binm{n+b}{n}\binm{n+c}{n}}
\sum_{i=0}^n(-1)^i\binm{n}{i}\frac{\binm{i+b}{i}\binm{i+c}{i}\binm{i+d+e+n}{i}}
{\binm{i+d}{i}\binm{i+e}{i}\binm{i+b+c+1}{i}}
\sum_{j=0}^i\binm{i}{j}\frac{\binm{n+d}{j}\binm{n+e}{j}\binm{j+f+g+n}{j}}
{\binm{j+f}{j}\binm{j+g}{j}\binm{i+d+e+n}{j}}.
  \enm
\end{exam}

\begin{exam}\label{exam-n}
 For six nonnegative integers $\{b,c,d,e,f,g\}$, there holds
\bnm
 &&\xqdn\! \sum_{k=0}^n\binm{n}{k}^2
 \frac{\binm{k+b}{k}\binm{n+c}{k}\binm{n+d}{k}\binm{n+e}{k}\binm{n+f}{k}\binm{n+g}{k}}
{\binm{n+b}{k}\binm{k+c}{k}\binm{k+d}{k}\binm{k+e}{k}\binm{k+f}{k}\binm{k+g}{k}}
\\&&\xqdn\!\:\:\times\:
\Big\{1+(n-2k)(2H_k-H_{b+k}+H_{c+k}+H_{d+k}+H_{e+k}+H_{f+k}+H_{g+k})\Big\}\\
&&\xqdn\!\:\:=\:(-1)^n\frac{\binm{b-c}{n}}{\binm{n+b}{n}\binm{n+c}{n}}
\sum_{i=0}^n\binm{n}{i}\frac{\binm{i+b}{i}\binm{n+c}{i}\binm{i+d+e+n}{i}}
{\binm{i+d}{i}\binm{i+e}{i}\binm{i+b-c-n}{i}}
\sum_{j=0}^i\binm{i}{j}\frac{\binm{n+d}{j}\binm{n+e}{j}\binm{j+f+g+n}{j}}
{\binm{j+f}{j}\binm{j+g}{j}\binm{i+d+e+n}{j}}.
  \enm
\end{exam}

\begin{exam}\label{exam-o}
 For six nonnegative integers $\{b,c,d,e,f,g\}$,  there holds
\bnm
 &&\xqdn \sum_{k=0}^n\binm{n}{k}^2
 \frac{\binm{n+b}{k}\binm{n+c}{k}\binm{n+d}{k}\binm{n+e}{k}\binm{n+f}{k}\binm{n+g}{k}}
{\binm{k+b}{k}\binm{k+c}{k}\binm{k+d}{k}\binm{k+e}{k}\binm{k+f}{k}\binm{k+g}{k}}
\\&&\xqdn\:\:\times\:
\Big\{1+(n-2k)(2H_k+H_{b+k}+H_{c+k}+H_{d+k}+H_{e+k}+H_{f+k}+H_{g+k})\Big\}\\
&&\xqdn\:\:=\:(-1)^n\frac{\binm{2n+b+c}{n}}{\binm{n+b}{n}\binm{n+c}{n}}
\sum_{i=0}^n\binm{n}{i}\frac{\binm{n+b}{i}\binm{n+c}{i}\binm{i+d+e+n}{i}}
{\binm{i+d}{i}\binm{i+e}{i}\binm{2n+b+c}{i}}
\sum_{j=0}^i\binm{i}{j}\frac{\binm{n+d}{j}\binm{n+e}{j}\binm{j+f+g+n}{j}}
{\binm{j+f}{j}\binm{j+g}{j}\binm{i+d+e+n}{j}}.
  \enm
\end{exam}

It should be pointed out that Examples \ref{exam-a}-\ref{exam-h} are
only the suitable limiting cases of Examples
\ref{exam-i}-\ref{exam-o}.
 Although the latter are also crossed one another as the former, they can create
numerous beautiful harmonic number identities with doubt. Further,
Theorem \ref{thm} may produce more harmonic number identities of
Chu-Donno type with the change of $m$. The interested reader may
write several ones of them down as exercises.

\subsection{Special cases: harmonic number identities of
Paule-Schneider type} ${}$

For an integer $u$ with $u\neq 0$, define $T_n^{(u)}$ by
\[T_n^{(u)}=\sum_{k=0}^n\binm{n}{k}^{u}\big\{1+u(n-2k)H_k\big\}.\]
Then eight known harmonic number identities can be stated as
follows:
 \bmn
T_n^{(-2)}&=&2\frac{(1+n)^2}{(2+n)}H_{n+1},        \label{har-a}\\
T_n^{(-1)}&=&(1+n)H_{n+1},                          \label{har-b}\\
T_n^{(1)}&=&1,                                        \label{har-c}\\
T_n^{(2)}&=&0,                                        \label{har-d}\\
T_n^{(3)}&=&(-1)^n,                                   \label{har-e}\\
T_n^{(4)}&=&(-1)^n\sbnm{2n}{n},                         \label{har-f}\\
T_n^{(5)}&=&(-1)^n\sum_{i=0}^n\sbnm{n}{i}^2\sbnm{n+i}{n}, \label{har-g}\\
T_n^{(6)}&=&(-1)^n\sum_{i=0}^n\sbnm{n}{i}^2\sbnm{n+i}{n}\sbnm{2n-i}{n}.
\label{har-h}
 \emn
\eqref{har-c}-\eqref{har-g} appeared  first in Paule and Schneider
\cito{paule}. Chu and Donno \cito{chu} offered other three ones and
showed that these eight harmonic number identities just displayed
can be derived by specifying the parameters in Examples
\ref{exam-a}, \ref{exam-c}, \ref{exam-d} and \ref{exam-h}. Now, we
shall bend ourselves to display the remaining results of the same
type by specifying the parameters in Theorem \ref{thm}.

Letting $P_{2m+2}\to\infty$, $v\to0$ and $P_s\to0$ with $1\leq
s\leq2m+1$ in Theorem \ref{thm}, we obtain the equivalent form of
the first equation of Krattenthaler and Rivoal
\citu{krattenthaler-a}{Proposition 1}.

\begin{prop}\label{prop-c}
 For $m\in\mathbb{N}$, there holds the harmonic number identity:
\bnm
 T_n^{(2m+3)}=(-1)^n\sum_{0\leq i_1\leq i_2\leq\cdots\leq i_m \leq n}
\binm{n}{i_m}^2\binm{n+i_1}{n}
 \prod_{r=1}^{m-1}\binm{n}{i_r}^2\binm{n+i_{r+1}-i_r}{n}.
 \enm
\end{prop}

Setting $v=0$ and $P_s=0$ with $1\leq s\leq2m+2$ in Theorem
\ref{thm}, we get the following equation.

\begin{prop}\label{prop-d}
 For $i_{m+1}=n$ with $m\in\mathbb{N}$, there holds the harmonic number identity:
\bnm
  T_n^{(2m+4)}=(-1)^n\sum_{0\leq i_1\leq i_2\leq\cdots\leq i_m \leq n}
\binm{n+i_1}{n}\prod_{r=1}^{m}\binm{n}{i_r}^2\binm{n+i_{r+1}-i_r}{n}.
 \enm
\end{prop}

Proposition \ref{prop-d} and the second equation of Krattenthaler
and Rivoal \citu{krattenthaler-a}{Proposition 1} have different
versions although that the purpose of them are the same. Proposition
\ref{prop-d} reduces to \eqref{har-h} exactly when $m=1$. Other two
results are laid out as follows.

\begin{exam}[Harmonic number identity of Paule-Schneider type: $m=2$ in Proposition \ref{prop-d}]
\bnm \xxqdn\xxqdn\qqdn
 T_n^{(8)} =(-1)^n\sum_{i=0}^n
\binm{n}{i}^2\binm{2n-i}{n}
\sum_{j=0}^i\binm{n}{j}^2\binm{n+j}{n}\binm{n+i-j}{n}.
 \enm
\end{exam}

\begin{exam}[Harmonic number identity of Paule-Schneider type: $m=3$ in Proposition \ref{prop-d}]
\bnm
 \quad\:\; T_n^{(10)}=(-1)^n\sum_{i=0}^n \binm{n}{i}^2\binm{2n-i}{n}
\sum_{j=0}^i\binm{n}{j}^2\binm{n+i-j}{n}
\sum_{t=0}^j\binm{n}{t}^2\binm{n+t}{n}\binm{n+j-t}{n}.
 \enm
\end{exam}

The open problem posed at the end of Paule and Schneider
\cito{paule} states that whether $ T_n^{(u)} $ can be expressed as a
definite hypergeometric single-sum for all $u\geq3$. Although the
equations on $T_n^{(u)}$ with $u\geq3$ have been given in
Krattenthaler and Rivoal \citu{krattenthaler-a}{Proposition 1} and
this subsection, we can't still judge that whether $T_n^{(u)}$ can
be expressed as a definite hypergeometric single-sum for all
$u\geq6$.

 Letting $v\to 2m+2$, $P_{2m+2}\to\infty$ and $P_s\to0$ with
$1\leq s\leq2m+1$ in Theorem \ref{thm}, we attain the following
equation.

\begin{prop}\label{prop-e}
 For $m\in\mathbb{N}$, there holds the harmonic number identity:
\bnm
 T_n^{(1-2m)}=(1+n)^m\sum_{0\leq i_1\leq i_2\leq\cdots\leq i_m \leq n}
\frac{1}{1+n-i_1}
 \prod_{r=1}^{m-1}\frac{(1)_{i_r}(-i_{r+1})_{i_r}}
 {(-n)_{i_r}(1+n-i_{r+1})_{i_r+1}}.
 \enm
\end{prop}

Proposition \ref{prop-e} leads to \eqref{har-b} exactly when $m=1$.
Other two results are displayed as follows.

\begin{exam}[Harmonic number identity of Paule-Schneider type: $m=2$ in Proposition \ref{prop-e}]
\bnm
 \xxqdn\xxqdn\xqdn\qdn\: T_n^{(-3)}
=(1+n)^2\sum_{i=0}^n\sum_{j=0}^i\frac{(1)_j(-i)_j}{(-n)_j(1+n-i)_{j+1}}
\frac{1}{1+n-j}.
 \enm
\end{exam}

\begin{exam}[Harmonic number identity of Paule-Schneider type: $m=3$ in Proposition \ref{prop-e}]
\bnm
   \qquad\:\: T_n^{(-5)}
=(1+n)^3\sum_{i=0}^n\sum_{j=0}^i\frac{(1)_j(-i)_j}{(-n)_j(1+n-i)_{j+1}}
\sum_{t=0}^j\frac{(1)_t(-j)_t}{(-n)_t(1+n-j)_{t+1}}\frac{1}{1+n-t}.
 \enm
\end{exam}

Taking $v=2m+2$ and $P_s=0$ with $1\leq s\leq2m+2$ in Theorem
\ref{thm}, we achieve the following equation.

\begin{prop}\label{prop-f}
 For $i_{m+1}=n$ with $m\in\mathbb{N}$,
 there holds the harmonic number identity:
\bnm
 \quad T_n^{(-2m)}=(1+n)^{m+1}\sum_{0\leq i_1\leq
i_2\leq\cdots\leq i_m \leq n}
 \prod_{r=1}^{m}\frac{1}
 {1+n-i_r}\fnk{cccc}{1,&-i_{r+1}}{-n,&2+n-i_{r+1}}_{i_r}.
 \enm
\end{prop}

Proposition \ref{prop-f} reduces to \eqref{har-a} exactly when
$m=1$. Other two results are laid out as follows.

\begin{exam}[Harmonic number identity of Paule-Schneider type: $m=2$ in Proposition \ref{prop-f}]
\bnm
 \xxqdn \xxqdn \xxqdn T_n^{(-4)}
=(1+n)^3\sum_{i=0}^n\frac{1}{1+i}\sum_{j=0}^i\frac{(1)_j(-i)_j}{(-n)_j(1+n-i)_{j+1}}
\frac{1}{1+n-j}.
 \enm
\end{exam}

\begin{exam}[Harmonic number identity of Paule-Schneider type: $m=3$ in Proposition \ref{prop-f}]
\bnm
   \quad T_n^{(-6)}
=(1+n)^4\sum_{i=0}^n\frac{1}{1+i}\sum_{j=0}^i\frac{(1)_j(-i)_j}{(-n)_j(1+n-i)_{j+1}}
\sum_{t=0}^j\frac{(1)_t(-j)_t}{(-n)_t(1+n-j)_{t+1}}\frac{1}{1+n-t}.
 \enm
\end{exam}



\end{document}